\documentclass[12pt]{article}
\usepackage{amsmath,amssymb,amscd}
\newtheorem{theorem}{Theorem} 
\newtheorem{lemma}{Lemma}
\newtheorem{cor}{Corollary}

\newcommand{\qed}{{\unskip\nobreak\hfil\penalty50\hskip2em\vadjust{}
\nobreak\hfil$\Box$\parfillskip=0pt\finalhyphendemerits=0\par}}

\newcommand{\ds}{\displaystyle}

\newcommand{\m}{multi}
\renewcommand{\o}{polynomial}

\newcommand{\bs}{\bigskip}
\renewcommand{\d}{\delta}
\newcommand{\g}{\gamma}
\newcommand{\p}{\partial}
\renewcommand{\th}{\theta}
\newcommand{\ef}{entire functions }

\begin{document}
\setlength{\baselineskip}{18pt}

\begin{center}{\large\bf  Gauss-Lucas Theorems for Entire Functions
on $C^M$}\\ \bigskip
{\bf Marek Kanter}\end{center}

\bs\bs
\begin{quote} {\bf Abstract.} A  Gauss-Lucas theorem is proved for
\m variate entire functions, using a natural notion of separate convexity to
obtain sharp results. Previous work in this
area is mostly restricted to univariate entire functions (of genus no greater
than one unless ``realness" assumptions are made). The present work applies
to multivariate \ef whose sections can be written as a monomial times a canonical product
of arbitrary genus. Essential use is  made of the Levy-Steinitz theorem
for conditionally convergent vector series, a result generalizing
Riemann's well known theorem for conditionally convergent real number series.

\vspace{.33truein}\noindent
entire function, multivariate function, stable polynomial,\\ Levy-Steinitz theorem,
conditionally convergent series.
\end{quote}

\vspace{.25truein}

 Secondary \ 12D10, \ 26B25, \ 26C10

\setlength{\baselineskip}{18pt}

\section{Introduction}

Let $f (z)$ be a non-constant univariate \o \ and let $f_{,1}(z)$
stand for the complex derivative $\frac{d}{d z}\, f(z)$. 
The classical Gauss-Lucas theorem is the set relation
\begin{equation}\label{1}
f^{(-1)}_{,1}(0) \subset H(f^{(-1)}(0)) \ ,
\end{equation}
where $f^{(-1)}_{,1}(0)\equiv \{z: f_{,1}(z)\!=\!0\}$ and
 $H(f^{(-1)}(0))$ is the convex hull of the roots of $f (z)$ in the
complex plane $C$. It is of interest to extend this elegant result
to entire functions on $C^M$.

The desired Gauss-Lucas relation for non-constant entire univariate
functions $f (z)$ is the set inclusion
\begin{equation}\label{2}
f^{(-1)}_{,1}(0) \subset {\bar H}(f^{(-1)}(0)) \ ,
\end{equation}
where $\bar H(f^{(-1)}(0))$ is the closure of $H(f^{(-1)}(0))$ in the
complex plane. (Closure is taken since $f^{(-1)}(0)$ is not
necessarily finite.) Theorem 1 (in Section 2) proves this relation
under three assumptions:

a) that $f (z)=z^q g(z)$, where $q$ is a non-negative integer and
$g(z)$ is a canonical product of finite genus $p$ with $g(0)\!=\!1$.

b) that $g(z)$ has only finitely many zeros in any bounded
subset of $C$.

c) that the non-zero roots $(\gamma_n: n\geq 1)$ of $g(z)$ satisfy
\begin{equation}\label{3}
\sum^\infty_{x_n\geq 0} x_n = \infty \ ,
\end{equation}
where $(x_n: n\geq 1)$ is either the sequence
$Re(\gamma^{-1}_n)$, or the sequence $Re(\gamma_n^{-2}),\dots$, or the
sequence $Re(\gamma_n^{-p})$, or any one of the remaining sequences
$Re(-\gamma_n^{-r})$, $Im(\gamma_n^{-r})$, $Im(-\gamma_n^{-r})$
for  $1\leq r\leq p$.

 Previous work extending the classical Gauss-Lucas theorem has been 
 limited to entire functions of genus no greater than one, unless
 ``realness" assumptions are made. (See e.g.~[2], [3], [7], or [8].)
 Here, these restrictions are removed by appealing to the Levy-Steinitz
 theorem for vector series, which generalizes Riemann's theorem 
 that a conditionally convergent series of real numbers
 (convergent but not absolutely convergent) can be 
 rearranged to sum to any real number.

To formulate the Gauss-Lucas relation for {\it multivariate} 
\ef \newline $f (z)=f (z_1,\dots ,z_M)$, the notation 
$f_{,m}(z)\equiv \frac{\p}{\p z_m} \,f (z)$ and
\[
f (w,z^{(m)}) \equiv f (z_1,\dots ,z_{m-1},w,z_{m+1},\dots ,z_M)
\]
is convenient, where $z^{(m)}\equiv (z_1,\dots ,z_{m-1},z_{m+1},\dots ,z_M)$
is in $C^{M-1}$. It is then possible to define $f (z)$ as being an
entire function on $C^M$ if $g(w)\equiv f (w,z^{(m)})$ is an entire 
function on $C$ for all $z^{(m)}\in C^{M-1}$. (See Section 3.)

The obvious notion of convexity in $C^M$, inherited by identifying
$C^M$ with $R^{2M}$, can be extended in a natural way to yield
sharper results at no cost to simplicity. Defining 
the projection operator $P_m(z)=z_m$,
a set $K\subset C^M$ is called {\it separately convex in} $C^M$
if the section
\[
K(z^{(m)}) \equiv P_m(\{ y: y\in K, \ y^{(m)}\!=\!z^{(m)}\})
\]
is a convex subset of $C$ for all $m\in \{1,\dots ,M\}$ and all
$z^{(m)}\in C^{M-1}$. The class of all separately convex subsets
of $C^M$ is clearly closed under intersection. Thus the
convex hull of a set $A\subset C^M$ with respect to this notion
of convexity may be identified as the smallest set containing $A$
and separately convex in $C^M$. It is denoted $H_2(A)$.

Assume the partial derivative $f_{,m}(z)$ of $f (z)$ is not
identically zero in $z$. The desired 
Gauss-Lucas relation for multivariate \ef $f (z)$ is the set
inclusion
\begin{equation}\label{4}
f^{(-1)}_{,m}(0) \subset \bar H_2(f^{(-1)}(0)) \ ,
\end{equation}
where $\bar H_2(f^{(-1)}(0))$ is the closure of $H_2(f^{(-1)}(0))$
in $C^M$. Theorem 2 (in Section 3) establishes (\ref{4})
for a general class of multivariate \ef of finite genus.
(A multivariate entire function $f (z)$ is said to have
finite genus if the univariate function 
$g(w)\equiv f (w,z^{(m)})$ has finite genus for all
$m\in\{1,\dots ,M\}$ and all $z^{(m)}\in C^{M-1}$.)

Any subset of $C^M$ which is convex in the usual sense
is also separately convex. Letting
 $H(f^{(-1)}(0))$ be the usual convex hull of
 $f^{(-1)}(0)$, it follows that
\begin{equation}\label{5}
 H_2(f^{(-1)}(0)) \subset H(f^{(-1)}(0)) \ ,
\end{equation}
 whence the set relation (\ref{4})
 is sharper than the corresponding relation using the
 standard notion of convexity in $C^M$. It is shown in Section 3
 that the sharper relation (\ref{4}) implies related recent
 results about multivariate stable polynomials. (See [10]
 for an exposition of these results, wherein a polynomial
 on $C^M$ is termed stable if it has no roots $z$ with
 $Im(z_m)>0$ for all $m$.) If $f$ is a multivariate polynomial
 then closure need not be taken in (\ref{4}), yielding
\begin{equation}\label{6}
f^{(-1)}_{,m}(0) \subset H_2(f^{(-1)}(0)) \ .
\end{equation}  
(See [5,\,Theorem 1] for this result and further comments
on separately convex subsets of~$C^M$.)

\section{Univariate entire functions}

This section relates the Gauss-Lucas theorem for univariate
\ef to the Levy-Steinitz generalization of Riemann's theorem for conditionally convergent real series.

An entire function $f(z)$ (with non-zero roots 
$\gamma=(\gamma_n: n\geq 1)$ and a root at the origin of multiplicity $q$)
is a canonical product if $f (z)=z^q g(z)$, where 
\[
g(z) = \prod^{\infty}_{n=1} \left(1- \left(\frac{z}{\gamma_n}\right)\!\!\right)
\exp (w_n(z;\gamma))
\]
for \ $w_n(z;\gamma)\equiv \ds{\sum^p_{r=1} r^{-1}\left( \frac{z}{\gamma_n}\right)^{\!\!r}}$.
The genus of $f (z)$ is the smallest non-negative integer $p$ such that the sum
\begin{equation}\label{7}
\sum^\infty_{n=1} |\gamma_n|^{-(1+p)}   
\end{equation}
is finite. Let
\[
V_N(r;\gamma)=\sum^N_{n=1} \gamma_n^{-r} \ .
\]
It will be shown that $f (z)$ satisfies the Gauss-Lucas relation
(\ref{2}) if the simple condition
\begin{equation}\label{8}
\lim_{N\to\infty} V_N(r;\gamma) =0
\end{equation}
holds for $r\in \{1,\dots ,p\}$. Noting that the left-hand side of (\ref{8})
represents an infinite sum of vectors in $C^p$ 
(with coordinates indexed by $r$) and that this sum may depend on the ordering
of the non-zero roots $\gamma_n$ of the canonical product
$g(z)$, it is natural to ask if there is some reordering of these
roots that makes (\ref{8}) hold.

The Levy-Steinitz theorem for conditionally convergent series of
vectors is relevant to this question. 
The treatment of the Levy-Steinitz theorem that yields (\ref{8})
is due to Katznelson and McGehee [6].
These authors consider vectors in $R^{\infty}$, so it is useful
to write $\gamma^{-r}_n$ in the form
$x_n(r;\gamma) +ix_n(r\!+\!p;\gamma)$, for $1\leq r\leq p$, where
$x_n(k;\gamma)\in R$, for $1\leq k\leq 2p$ \ and \ $i=\sqrt{-1}$.

\begin{lemma} Suppose $\g\!=\!(\g_n: n\geq 1)$ is a sequence of
complex numbers such that $(\ref{3})$ holds. Then there exists a
permutation $\pi$ of the integers $(n\geq 1)$ such that
the rearrangement $\delta_n\equiv \g_{\pi (n)}$ satisfies $(\ref{8})$
if $\delta\!=\!(\delta_n: n\geq 1)$ is substituted for $\g$.
(Note the roots stay the same, only their ordering changes.)
\end{lemma}

{\bf Proof.} It follows from Riemann's theorem that for each 
$k\in \{1,\dots ,2p\}$, the series
\[
X(k\,;\g)\equiv \sum^\infty_{n=1} x_n(k\,;\g) 
\]
can  be rearranged to sum to any real number, {\it in particular
to the real number} 0. (The rearrangement may depend on $k$.)
 Let $S_X$ stand for the subset of
$R^{2p}$ defined by $\{(X(k;\delta): 1\!\leq \!k\!\leq\! 2p), \delta\!\in\!\Lambda\}$, where
$\Lambda$ consists of all rearrangements 
$\d\!=\!(\d_n: n\geq 1)$   of $\gamma =(\g_n: n\!\geq\! 1)$  that yield a
convergent sum for the vector series $\sum^\infty_{n=1} X_n(\d)$
in $R^{2p}$ with summands
\[
X_n(\d)\equiv (x_n(k\,;\d): 1\leq k\leq 2p) \ .
\]
It follows from Theorem 1 of Katznelson and McGehee that $S_X$ is
the solution of a finite number of homogeneous linear equations
in $R^{2p}$. In particular, $S_X$ is  a linear subspace of $R^{2p}$
and contains the zero vector. Translating back from $R^{2p}$ to $C^p$,
it follows that 
(\ref{8}) is satisfied if $\g$ is replaced by any rearrangement  $\d$
corresponding to the zero vector in $S_X$. \qed

\bs{\bf Remark 1.} The hypothesis of Lemma 1 can be replaced with
the statement that the sequence $(\g^{-r}_n: n\geq 1)$ has a
rearrangement summing to 0, for each $r\in \{1,\dots ,p\}$. This
modification allows the cases when $Re(\g^{-r}_n)$ or $Im(\g^{-r}_n)$
is absolutely summable to 0. Note further that the original
hypothesis is consistent with this modification, as follows from
Riemann's theorem and the work of Katznelson and 
McGehee. (Without the latter work, it is not clear that a sequence
$x_n+iy_n$ is rearrangeably summable to zero
if the same holds separately for $x_n$ and $y_n$.)

\bs{\bf Definition 1.} Let  $f (z)=z^q g(z)$, where $g(z)$ is a canonical
product of genus $p$. Then  $f(z)$ is
{\it rearrangeable} if the (non-zero) roots 
of $g(z)$ can be rearranged
so as to satisfy (\ref{8}).

\begin{theorem} If the entire function $f (z)=z^q g(z)$ is
rearrangeable, then it satisfies the Gauss-Lucas relation 
$(\ref{2})$.
\end{theorem}

{\bf Proof.} To start, assume $q\!=\!0$. It is desired to show that
\begin{equation}\label{9}
g(z)=\lim_{N\to\infty} f_N(z;\d)
\end{equation}
uniformly on compact subsets of $C$, where
\[
f_N(z;\d) \equiv \prod^N_{n=1} (1-(z/\d_n)) 
\]
and  $\d\!=\!(\d_n: n\geq 1)$ is a rearrangement of  the roots
of $g(z)$ satisfying (\ref{8}).  Let
\[
h_N(z;\d) =\sum^p_{r=1} r^{-1} V_N(r;\d)z^r \ .
\]
It is shown in Ahlfors [1, p.193] that the partial canonical
product
\[
g_N(z;\d) =\left(\prod^N_{n=1} (1-(z/\d_n))\right)\exp(h_N(z;\d))
\]
is absolutely and uniformly convergent to $g(z)$ on compact
subsets of $C$. (Thus $g(z)$ is rearrangement-invariant,
i.e. it does not depend on the ordering of its roots.) Note
$h_N(z;\d)$ converges uniformly to 0 on compact subsets of $C$
because the rearrangement $\d$ satisfies (\ref{8}). Furthermore (\ref{9})
holds since $|f_N(z;\d)-g(z)|$ is dominated by
\[
|(\exp (-h_N(z;\d))-1)g_N(z)| + |g_N(z)-g(z)|
\]
(the triangle inequality for complex numbers).

It remains to show that $f(z)$ satisfies (\ref{2}). Assume,
by way of contradiction, that there exists $z_0\in C$ with
$f'(z_0)=0$ and $z_0\not\in {\bar H}(f^{(-1)}(0))$. Let $K$
be a closed convex subset of $C$ with $z_0\not\in K$ and
$f^{(-1)}\!(0)\subset K$. By a theorem due to Hurwitz there 
exists a sequence of complex numbers $v_n$ tending to $z_0$
and a sequence of positive integers $M_n$ tending to infinity
such that $f'_{M_n}\!(v_n)=0$ for all~$n$. (See Titchmarsh [9].)
It follows (by the Gauss-Lucas theorem for polynomials
and the fact that $f^{(-1)}_N(0) \subset f^{(-1)}(0)$ for all $N$)
that $v_n\in K$ for all $n$. This contradicts the convergence of
$v_n$ to $z_0$.

The restriction that $q\!=\!0$ can be removed because the
argument that $f(z)$ satisfies the Gauss-Lucas relation 
(\ref{2}) is not affected by multiplication by the monomial $z^q$. \qed

\bs {\bf Remark 2.} It should be noted that polynomials and all
other entire functions of genus zero are rearrangeable because
condition (\ref{8}) is vacuous if $p\!=\!0$.

\section{Multivariate \ef }

In this section it will be shown that the multivariate
Gauss-Lucas relation (\ref{4}) follows  by successively holding
constant all variables but one.  Hormander [4] is the standard
reference about multivariate \ef  (e.g. Hartog's theorem
that separate analyticity implies joint analyticity).

\begin{lemma} Suppose $f(z)=f(z_1,\dots ,z_M)$ is an entire
function such that for all $m\!\in \!\{1,\dots ,M\}$ and all
$z^{(m)}\in C^{M-1}$, the univariate entire function
$g(w)\equiv f (w,z^{(m)})$ satisfies the univariate
Gauss-Lucas relation $(\ref{2})$. Then $f(z)$ satisfies the
multivariate Gauss-Lucas relation $(\ref{4})$. 
\end{lemma}

{\bf Proof.}  Suppose $f_{,m}(z)=0$. For $g(w)$ as above,
note that the derivative
\[
g'(z_m) = f_{,m}(z_m,z^{(m)}) = f_{,m}(z) \ .
\]
It follows $g'(z_m)\!=\!0$ and thus $z_m\in {\bar H}(g^{(-1)}(0))$.
Let $K$ be any closed separately convex subset of $C^M$
containing $f^{(-1)}(0)$. The section $K(z^{(m)})$ is a closed 
convex subset of $C$ containing $g^{(-1)}(0)$. Thus
$K(z^{(m)})$ contains ${\bar H}(g^{(-1)}(0))$. In particular
$z_m\in K(z^{(m)})$, i.e. $z\in K$. By choice of $K$ this
shows $z\in {\bar H}_2(f^{(-1)}(0))$. \qed

\bs {\bf Remark 3.} If the hypothesis in Lemma 2 is restated
as ``$g(w)$ satisfies the Gauss-Lucas relation (\ref{1})", then
the conclusion becomes ``$f(z)$ satisfies the 
multivariate Gauss-Lucas relation (\ref{6})"; the  proof is
the same except that the closures of $H(g^{-1}(0))$ and
$H_2(f^{(-1)}(0))$ are not taken. This change reflects the
simple argument in Kanter [5] concerning multivariate
polynomials.

\begin{theorem} Suppose $f(z_1,\dots ,z_M)$ is an entire 
function such that the univariate function
$g(w)=f(w, z^{(m)})$ is rearrangeable for all $m\!\in\! \{1,\dots ,M\}$
and all $z^{(m)}\in C^{M-1}$. Then $f(z)$ satisfies
the Gauss-Lucas relation $(\ref{4})$. 
\end{theorem}

{\bf Proof.} Apply Theorem 1 and Lemma 2. \qed

\bs Theorem 2 can be immediately applied to multivariate stable 
\ef \!\!.

\bs{\bf Definition 2.} Let $\th =(\th_1,\dots \th_M)\in R^M$
and let 
\[
A(\th) =\{ z\in C^M: {\rm{Im}}(e^{i\th_m} z_m) >0 \  {\mbox{for }} 1\leq m\leq M\} \ .
\]
A multivariate entire function is called $\th$-stable if it has
no zeros in $A(\th)$.

\bs{\bf Remark 4.} Let $A^c(\th)$ stand for the complement
of $A(\th)$ in $C^M$. It is easy to see that $A^c(\th)$
is a closed separately convex subset of $C^M$. (See [5].)

\begin{cor} Suppose the multivariate $\th$-stable entire
function $f(z)$ satisfies the hypotheses of Theorem $2$.
Then any non-null partial derivative $f_{,m}(z)$ is also
$\th$-stable.
\end{cor}

\newpage
\begin{center}{\bf References}\end{center}

\begin{description}  
\item{[1]} \ L.\! V.\! Ahlfors, {\it Complex Analysis}, McGraw-Hill,
New York, 1966.

\item{[2]} \ J.\! Clunie and Q.\! I.\! Rahman, Extension of a
theorem of J.H.Grace to transcendental \ef \!\!, {\it
Mathematical Proceedings of the Cambridge Philosophical
Society} {\bf 112} (1992), 565--573.

\item{[3]} \ T.\! Craven and G.\! Csordas, The Gauss-Lucas theorem
and Jensen polynomials, {\it Trans. Amer. Math. Soc.} {\bf 357}
(2004), 4043--4063.

\item{[4]} \ L.\! Hormander, {\it An Introduction to Complex Analysis
in Several Variables}, North Holland Publishing Co., New York, 1973.

\item{[5]} \ M.\! Kanter, Multivariate Gauss-Lucas theorems,\\
http:/\!/arxiv.org/1203.6426v1 .

\item{[6]} \ V.\! Katznelson and O.\! C.\! McGehee, Conditionally
convergent series in $R^{\infty}$, {\it Michigan Math. Journal}
{\bf 21} (1974), 97--106.

\item{[7]} \ M.\! Marden, On the zeros of the derivative of an
entire function, {\it The American Math. Monthly} {\bf 75} 
(1968), 829--839.

\item{[8]} \ M.\! B.\! Porter, On a theorem of Lucas, {\it Proc.
Nat. Acad. Sci. USA} {\bf 2} (1916), 247--248, 335--336.

\item{[9]} \ E.\! C.\! Titchmarsh, {\it The Theory of Functions},
Oxford University Press, London, 1979.

\item{[10]} \ D.\! G.\! Wagner, Multivariate stable polynomials:
theory and applications, {\it Bulletin Amer. Math. Soc.} {\bf 48}
(2011), 53--83.
\end{description}

\bs\bs\hfill{}\begin{tabular}{r}
{\sc 1216 Monterey Avenue, Berkeley, California 94707}\\
{\it E-mail address}  \ {\tt mrk$@$cpuc.ca.gov}
\end{tabular}

\end{document}